\documentclass[a4paper,reqno,12pt]{amsart}
\usepackage[utf8]{inputenc}
\usepackage[T1,T2A]{fontenc}
\usepackage[english]{babel}
\usepackage{amsmath}
\usepackage{amssymb}
\usepackage{amscd}
\usepackage{amsthm}
\usepackage{euscript}
\usepackage{enumitem}
\usepackage[all]{xy}
\newtheorem{proposition}{Proposition}[section]
\newtheorem{conjecture}[proposition]{Conjecture}

\newtheorem{question}[proposition]{Question}
\newtheorem{lemma}[proposition]{Lemma}
\newtheorem{theorem}[proposition]{Theorem}

\newtheorem{corollary}[proposition]{Corollary}
\theoremstyle{definition}
\newtheorem{definition}[proposition]{Definition}
\newtheorem{example}[proposition]{Example}
\theoremstyle{remark}

\DeclareMathOperator{\ad}{ad}

\DeclareMathOperator{\Der}{Der}
\DeclareMathOperator{\Char}{char}
\renewcommand{\Im}{\mathop{\mathrm{Im}}}
\def\FF{{\mathbb F}}

\def\RRR{\mathfrak{R}}
\def\SSS{\mathfrak{S}}

\def\ggg{\mathfrak{g}}
\def\hhh{\mathfrak{h}}
\def\FFF{\mathcal{F}}

\def\sl{\mathfrak{sl}}
\def\sp{\mathfrak{sp}}
\def\sch{\mathfrak{sch}}
\renewcommand{\phi}{\varphi}

\renewcommand{\le}{\leqslant}
\renewcommand{\phi}{\varphi}

\renewcommand{\le}{\leqslant}

\begin{document}
	\date{}
	\title[Uniqueness of addition in {Schr{\"o}dinger} algebras]{Uniqueness of addition in {Schr{\"o}dinger} algebras}
	\author{Gennadiy Sosnov}
	\address{Faculty of Computer Science, HSE University, Pokrovsky Boulevard 11, Moscow, 109028 Russia}
	\email{grsosnov@edu.hse.ru}
	\thanks{This article is an output of a research project (HSE-BR-2025-22) implemented as part of the Basic Research Program at HSE University.}
    
	\begin{abstract}
    We obtain a criterion for the uniqueness of addition in Schr{\"o}dinger Lie algebras, applicable to both bijections and injections. This allows us to construct the first example of a Lie algebra with unique addition and a non-trivial center over a field of characteristic zero. Furthermore, a necessary and sufficient condition is established for the Schr{\"o}dinger algebras to have bracket width one.
	\end{abstract}
	\maketitle


	\section{Introduction}\label{sec:1}

	A Lie ring $\RRR$ is called a Lie ring with unique addition, or a UA-Lie ring, if every Lie ring bijection $\phi: \RRR \to \SSS$ that preserves the commutator is additive, that is, for all $x, y \in \RRR$, we have $\phi(x + y) = \phi(x) + \phi(y)$.

	In~\cite[Proposition~1]{Arzhantsev2025} it is proved that if a Lie ring $\RRR$ has more than $4$ elements, $[\RRR, \RRR] \neq \RRR$, $Z(\RRR) \neq 0$, and the condition $Z(\RRR) \cap [\RRR, \RRR] = 0$ implies $|Z(\RRR)| > 2$, then $\RRR$ is not a UA-Lie ring. It is natural to ask whether the uniqueness of addition fails if these conditions are relaxed. The results from~\cite{Arzhantsev2000,Arzhantsev2001,Mayorova2019} do not suffice to answer this question. In \cite[Problem~3]{Arzhantsev2025}, the question of the uniqueness of addition in the Schr{\"o}dinger algebra $\mathfrak{sch}_1(K)$ is raised. This is a Lie algebra of minimal dimension that is perfect, meaning it coincides with its derived subalgebra, and possesses a non-trivial center.

	Let us recall that the Schr{\"o}dinger algebra is the Lie algebra of the Schr{\"o}dinger group, that is, the symmetry group of the configuration space of the Schr{\"o}dinger equation~\cite{Perroud1977}. In this article, we consider the $n$-th Schr{\"o}dinger algebra $\sch_n(K) = \sl_2(K) \ltimes \hhh_n(K)$, where $\hhh_n(K)$ is a $(2n + 1)$-dimensional Heisenberg algebra. These algebras have been actively studied recently. In~\cite{Tao2022}, the author obtains a classification of simple Harish-Chandra modules and simple Whittaker modules over the $n$-th Schr{\"o}dinger algebra. In~\cite{Zhengxin2023}, all derivations and biderivations of the algebra $\sch_n(K)$ are determined. In~\cite{Boqing2024}, the Lie algebra of derivations and the automorphism group of the algebra $\sch_n(K)$ are described.

	Replacing bijections by injections in the definition of a UA-Lie ring, one obtains the notion of an iUA-Lie ring. It turns out that any iUA-Lie ring is perfect and has bracket width one, that is, any element of the ring can be expressed as a single commutator~\cite[Lemma~4.3]{Sosnov2026}. 
  
  We will denote by $\FF_q$ the finite field with $q$ elements. Let us state the main result of the paper.

	\begin{theorem}\label{theorem:when_schrodinger_UA}
    Let $K$ be a field and $n$ be a positive integer.
    \begin{enumerate}[label=(\arabic*)]
      \item The Lie algebra $\sch_n(K)$ is a UA-Lie ring if and only if it is an iUA-Lie ring or~$\sch_1(\FF_3)$.
      \item If $\Char K = 2$, then the Lie algebra $\sch_n(K)$ is not a UA-Lie ring.
      \item If $\Char K \neq 2$, then the following statements are equivalent:
      \begin{enumerate}[label=(\alph*)]
        \item\label{item:1} Every element of $K$ is a sum of $n$ squares;
        \item\label{item:2} The algebra $\sch_n(K)$ has bracket width one;
        \item\label{item:3} The algebra $\sch_n(K)$ is an iUA-Lie ring.
      \end{enumerate}
    \end{enumerate}
	\end{theorem}
	
	The author is grateful to Ivan Arzhantsev for his valuable comments and assistance in the preparation of this work.

	\section{Schr{\"o}dinger algebra}\label{sec:2}

	Let us begin with the main definitions.

	\begin{definition}
		The \textit{Heisenberg algebra} $\hhh_n(K)$ over the field $K$ is a Lie algebra consisting of matrices of the form
		\begin{equation*}
      \begin{pmatrix}
        0 & \mathbf{a} & c \\
        \mathbf{0} & 0_n & \mathbf{b} \\
        0 & \mathbf{0} & 0
      \end{pmatrix},
		\end{equation*}
		where $0_n$ is the $n \times n$ zero matrix, $\mathbf{a}^{\top}$ and $\mathbf{b}$ are column vectors in $K^n$, and $c \in K$.
	\end{definition}

	It follows from the definition that the Heisenberg algebra has the basis
	\begin{equation*}
    x_i = 
    \begin{pmatrix}
      0 & e_i^{\top} & 0 \\
      \mathbf{0} & 0_n & \mathbf{0} \\
      0 & \mathbf{0} & 0
    \end{pmatrix}, \quad
    z = 
    \begin{pmatrix}
      0 & \mathbf{0} & 1 \\
      \mathbf{0} & 0_n & \mathbf{0} \\
      0 & \mathbf{0} & 0
    \end{pmatrix}, \quad
    y_j =
    \begin{pmatrix}
      0 & \mathbf{0} & 0 \\
      \mathbf{0} & 0_n & e_j \\
      0 & \mathbf{0} & 0  
    \end{pmatrix}
  \end{equation*}
	with relations
	\begin{equation*}
    [x_i, y_j] = \delta_{ij}z, \quad [x_i, z] = [y_j, z] = [x_i, x_j] = [y_i, y_j] = 0,
  \end{equation*}
	where $\delta_{ij}$ is the Kronecker delta and $e_i$ is an element of the standard basis of $K^n$.

	Let us fix the standard basis of $\sl_2(K)$:
	\begin{equation*}
    f =
    \begin{pmatrix}
      0 & 0 \\
      1 & 0
    \end{pmatrix}, \quad
    h = 
    \begin{pmatrix}
      1 & 0 \\
      0 & -1
    \end{pmatrix}, \quad
    e =
    \begin{pmatrix}
      0 & 1 \\
      0 & 0
    \end{pmatrix}.
  \end{equation*}

	\begin{definition}
		The \textit{Schr{\"o}dinger algebra} $\sch_n(K)$ is defined as the semidirect product $\sl_2(K) \ltimes \hhh_n(K)$, where the homomorphism $\sl_2(K) \rightarrow \Der(\hhh_n(K))$ is given by the relations
		\begin{equation*}
      \begin{matrix}
        [h, x_i] = x_i, & [h, y_j] = -y_j, & [h, z] = 0, \\
        [e, x_i] = 0, & [e, y_j] = x_j, & [e, z] = 0, \\
        [f, x_i] = y_i, & [f, y_j] = 0, & [f, z] = 0.
      \end{matrix}
    \end{equation*}
	\end{definition}

	The Schr{\"o}dinger algebra $\sch_n(K)$ is a perfect Lie algebra, that is, $[\sch_n(K), \sch_n(K)] = \sch_n(K)$ over a field $K$ of characteristic not $2$. The center of the algebra $\sch_n(K)$ is always equal to $\langle z \rangle$.

	The Schr{\"o}dinger algebra $\sch_n(K)$ admits the following matrix representation:
  \begin{equation*}
    f \mapsto
    \begin{pmatrix}
      0 & \mathbf{0} & \mathbf{0} & 0 \\
      \mathbf{0} & 0_n & 0_n & \mathbf{0} \\
      \mathbf{0} & I_n & 0_n & \mathbf{0} \\
      0 & \mathbf{0} & \mathbf{0} & 0
    \end{pmatrix},
    h \mapsto
    \begin{pmatrix}
      0 & \mathbf{0} & \mathbf{0} & 0 \\
      \mathbf{0} & I_n & 0_n & \mathbf{0} \\
      \mathbf{0} & 0 & -I_n & \mathbf{0} \\
      0 & \mathbf{0} & \mathbf{0} & 0
    \end{pmatrix},
    e \mapsto
    \begin{pmatrix}
      0 & \mathbf{0} & \mathbf{0} & 0 \\
      \mathbf{0} & 0_n & I_n & \mathbf{0} \\
      \mathbf{0} & 0_n & 0_n & \mathbf{0} \\
      0 & \mathbf{0} & \mathbf{0} & 0
    \end{pmatrix},
  \end{equation*}
  \begin{equation*}
    x_i \mapsto
    \begin{pmatrix}
      0 & \mathbf{0} & -e_i^{\top} & 0 \\
      \mathbf{0} & 0_n & 0_n & -e_i \\
      \mathbf{0} & 0_n & 0_n & \mathbf{0} \\
      0 & \mathbf{0} & \mathbf{0} & 0
    \end{pmatrix},
    z \mapsto
    \begin{pmatrix}
      0 & \mathbf{0} & \mathbf{0} & 2 \\
      \mathbf{0} & 0_n & 0_n & \mathbf{0} \\
      \mathbf{0} & 0_n & 0_n & \mathbf{0} \\
      0 & \mathbf{0} & \mathbf{0} & 0
    \end{pmatrix},
    y_j \mapsto
    \begin{pmatrix}
      0 & e_j^{\top} & \mathbf{0} & 0 \\
      \mathbf{0} & 0_n & 0_n & \mathbf{0} \\
      \mathbf{0} & 0_n & 0_n & -e_j \\
      0 & \mathbf{0} & \mathbf{0} & 0
    \end{pmatrix},
  \end{equation*}
	where $I_n$ is the identity matrix of order $n$. A direct verification shows that every matrix $X$ of this representation satisfies the relation $X^{\top}S + SX = 0$, where
	\begin{equation*}
		S =
		\begin{pmatrix}
			0 & 0 & 0 & 1 \\
			0 & 0_n & I_{n} & 0 \\
			0 & -I_{n} & 0_n & 0 \\
			-1 & 0 & 0 & 0
		\end{pmatrix},
	\end{equation*}
	that is, the Schr{\"o}dinger algebra $\sch_n(K)$ is embedded in the symplectic Lie algebra $\sp_{2n + 2}(K)$.

	\section{Bracket width}\label{sec:3}

  One of the characteristics of a Lie algebra $\ggg$ is the \textit{bracket width}, that is, the supremum of the values $\ell(a)$ over all $a$ in $[\ggg,\ggg]$, where $\ell(a)$ is defined as the minimal number $m$ in the decomposition of $a$ into a sum of commutators
  \begin{equation*}
    a = \sum_{i=1}^{m} [x_i,y_i].
  \end{equation*}

  It is known that the bracket width of finite-dimensional simple Lie algebras over an algebraically closed field of characteristic $0$ is equal to one~\cite{Brown1963}. Over an arbitrary field this question remains open; however, it is known that the bracket width does not exceed two. In~\cite{Dubouloz2021}, an example of an infinite-dimensional simple Lie algebra with bracket width strictly greater than one is constructed.

  We will be interested in the value of the bracket width of the perfect algebra $\sch_n(K)$. Since for $\Char K = 2$ the algebra $\sch_n(K)$ is not perfect, we will assume that $\Char K \neq 2$.

  \begin{proposition}\label{proposition:whenCommutatorWidthIs1}
    Let $K$ be a field with $\Char K \neq 2$. Then the bracket width of the algebra $\sch_n(K)$ does not exceed two and is equal to one if and only if every element $k \in K$ can be represented as a sum of $n$ squares.
  \end{proposition}

  To prove this, we will need some auxiliary lemmas.

  \begin{lemma}\label{lemma:commutator_width}
    Suppose that in the field $K$ every element can be represented as a sum of $n$ squares. Then the algebra $\sch_n(K)$ has bracket width one.
  \end{lemma}

  \begin{proof}
    Consider an arbitrary element
    \begin{equation*}
      C = Z + \sum_{i=1}^{n} \gamma_i x_i + \sum_{j=1}^{n} \gamma'_j y_j + cz, \quad \text{where} \quad Z =
      \begin{pmatrix}
        c_{11} & c_{12} \\
        c_{21} & c_{22}
      \end{pmatrix}.
    \end{equation*}
    We will look for a representation $C = [A,B]$, where
    \begin{equation*}
      A = X + \sum_{i=1}^{n} \alpha_i x_i + \sum_{j=1}^{n} \alpha'_j y_j \quad \text{and} \quad B = Y + \sum_{i=1}^{n} \beta_i x_i + \sum_{j=1}^{n} \beta'_j y_j.
    \end{equation*}

    The condition $C = [A,B]$ gives $n$ systems
    \begin{equation}\label{eq:systems}
      \begin{pmatrix}
        X & -Y
      \end{pmatrix}
      \begin{pmatrix}
        \beta_i \\ \beta'_i \\ \alpha_i \\ \alpha'_i
      \end{pmatrix} =
      \begin{pmatrix}
        \gamma_i \\
        \gamma'_i
      \end{pmatrix},
    \end{equation}
    the equation $[X,Y] = Z$, and the quadratic equation
    \begin{equation}\label{eq:quadratic}
      \sum_{i=1}^{n} \alpha_i \beta'_i - \beta_i \alpha'_i = c.
    \end{equation}

    \textit{Case 1.} Let $Z = 0$. If all $\gamma_i$ and $\gamma'_i$ are equal to $0$, then $A = cx_1$ and $B = y_1$. Otherwise, without loss of generality, we may assume that $\gamma_1 \neq 0$. Then $A = h - \frac{c}{\gamma_1}y_1$ and
    \begin{equation*}
      B = \gamma_1 x_1 - \gamma'_1 y_1 + \gamma_2 x_2 - \gamma'_2 y_2 + \cdots + \gamma_n x_n - \gamma'_n y_n.
    \end{equation*}

    \textit{Case 2.} Let $Z = \mu h$ and $\mu \neq 0$. Then take $X = \mu e$, $Y = f$. Systems~\eqref{eq:systems} give the solution $\beta'_i = \frac{\gamma_i}{\mu}$, $\alpha_i = -\gamma'_i$. Substituting into~\eqref{eq:quadratic} we obtain an equation which obviously has a solution.

    \textit{Case 3.} Let $c_{21} \neq 0$. Multiplying $C$ by $\frac{1}{c_{21}}$, we may assume that $c_{21} = 1$. Let
    \begin{equation*}
      X =
      \begin{pmatrix}
        1 & -2c_{11} \\
        0 & -1
      \end{pmatrix}, \quad
      Y =
      \begin{pmatrix}
        0 & \frac{c_{12}}{2} \\
        -\frac{1}{2} & 0
      \end{pmatrix}.
    \end{equation*}
    Systems~\eqref{eq:systems} give the following expressions
    \begin{equation*}
      \beta'_i = \frac{1}{2}\alpha_i - \gamma'_i, \quad \beta_i = \gamma_i - \frac{c_{12}}{2}\alpha'_i + \frac{2c_{11}}{c_{21}}\beta'_i = \gamma_i - \frac{c_{12}}{2}\alpha'_i + c_{11}\alpha_i - 2c_{11}\gamma'_i,
    \end{equation*}
    therefore equation~\eqref{eq:quadratic} takes the form
    \begin{equation*}
      \sum_{i=1}^{n} \frac{1}{2}\alpha_i^2 - \gamma'_i\alpha_i - \gamma_i\alpha'_i + \frac{c_{12}}{2}(\alpha'_i)^2 - c_{11}\alpha_i\alpha'_i + 2c_{11}\gamma'_i\alpha'_i = c.
    \end{equation*}
    Setting $\alpha'_i$ equal to one for all $i$, we obtain
    \begin{equation*}
      \sum_{i=1}^{n} \frac{1}{2}\alpha_i^2 - (\gamma'_i + c_{11})\alpha_i - \gamma_i + \frac{c_{12}}{2} + 2c_{11}\gamma'_i = c,
    \end{equation*}
    which is equivalent to
    \begin{multline*}
      \sum_{i=1}^{n} (\alpha_i - \gamma'_i - c_{11})^2 = \sum_{i=1}^{n} \alpha_i^2 - 2(\gamma_i' + c_{11})\alpha_i + (\gamma_i' + c_{11})^2 = \\
      = 2c + \sum_{i=1}^{n} 2\gamma_i - c_{12} - 4c_{11}\gamma'_i + (\gamma_i' + c_{11})^2,
    \end{multline*}
    hence one can choose $\alpha_i$ so that the equality holds.

    \textit{Case 4.} Let $c_{12} \neq 0$. Then one can use the previous case, since there is an involutive isomorphism of $\sl_2(K)$ given by the formula $A \mapsto -A^{\top}$, which extends to $\sch_n(K)$.
  \end{proof}

  \begin{lemma}\label{lemma:when_elem_is_not_commutator}
    Suppose that in the field $K$ the element $k$ cannot be represented as a sum of $n$ squares. Then the element $e - \frac{k}{2}z \in \sch_n(K)$ is not a commutator.
  \end{lemma}

  \begin{proof}
    Assume the contrary. Then there exist
    \begin{equation*}
      A = X + \sum_{i=1}^{n} \alpha_i x_i + \sum_{j=1}^{n} \alpha'_j y_j \quad \text{and} \quad B = Y + \sum_{i=1}^{n} \beta_i x_i + \sum_{j=1}^{n} \beta'_j y_j,
    \end{equation*}
    such that $[X,Y] = e$, for all $i$ we have
    \begin{equation*}
      \begin{pmatrix}
        X & -Y
      \end{pmatrix}
      \begin{pmatrix}
        \beta_i \\ \beta'_i \\ \alpha_i \\ \alpha'_i
      \end{pmatrix} =
      \begin{pmatrix}
        0 \\
        0
      \end{pmatrix} \quad \text{and} \quad \sum_{i=1}^{n} \alpha_i \beta'_i - \beta_i \alpha'_i = -\frac{k}{2}.
    \end{equation*}
    Then $X$ and $Y$ have the form
    \begin{equation*}
      X =
      \begin{pmatrix}
        a & b \\
        0 & -a
      \end{pmatrix}, \quad
      Y =
      \begin{pmatrix}
        a' & b' \\
        0 & -a'
      \end{pmatrix},
    \end{equation*}
    and moreover $2(ab' - a'b) = 1$. Then
    \begin{equation*}
      a'\alpha'_i - a\beta'_i = 0.
    \end{equation*}
    Since $a$ and $a'$ cannot both be zero, there exists an element $\lambda_i \in K$ such that $\alpha'_i = \lambda_i a$ and $\beta'_i = \lambda_i a'$. Substituting these values into the first equation of the system, we obtain
    \begin{equation*}
      a\beta_i - a'\alpha_i = \lambda_i(b'a - ba') = \frac{\lambda_i}{2}.
    \end{equation*}
    Now substitute $\alpha'_i$ and $\beta'_i$ into the expression with $-\frac{k}{2}$. We obtain
    \begin{equation*}
      \sum_{i=1}^{n} \lambda_i(a'\alpha_i - a\beta_i) = \sum_{i=1}^{n} -\frac{\lambda_i^2}{2} = -\frac{k}{2},
    \end{equation*}
    that is, $k$ is represented as a sum of $n$ squares, a contradiction.
  \end{proof}

  \begin{proof}[Proof of Proposition~\ref{proposition:whenCommutatorWidthIs1}]
    From Lemmas~\ref{lemma:commutator_width} and~\ref{lemma:when_elem_is_not_commutator} it follows that if the field $K$ contains an element which cannot be represented as a sum of $n$ squares, then the bracket width of $\sch_n(K)$ is greater than one. Therefore it suffices to prove that every element of $\sch_n(K)$ can be represented as a sum of two commutators.

    Consider an arbitrary element
    \begin{equation*}
      C = Z + \sum_{i=1}^{n} \gamma_i x_i + \sum_{j=1}^{n} \gamma'_j y_j + cz \in \sch_n(K).
    \end{equation*}
    From the proof of Lemma~\ref{lemma:commutator_width} it follows that there exist matrices $X, Y \in \sl_2(K)$ such that $Z = [X,Y]$ and the matrix $X$ is nondegenerate. Therefore
    \begin{equation*}
      C = [X, Y + \sum_{i=1}^{n} \alpha_i x_i + \alpha'_i y_i] + [cx_1, y_1],
    \end{equation*}
    where $\alpha_i$ and $\alpha'_i$ are chosen so that the identity
    \begin{equation*}
      X
      \begin{pmatrix}
        \alpha_i \\
        \alpha'_i
      \end{pmatrix} =
      \begin{pmatrix}
        \gamma_i \\
        \gamma'_i
      \end{pmatrix}
    \end{equation*}
    holds.
  \end{proof}

  We will also need a purely technical lemma.

  \begin{lemma}\label{lemma:commutatorOfSch1F3}
    Let $A = X + \lambda x_1 + \mu y_1 \in \sch_1(\FF_3)$, where $X \in \sl_2(\FF_3)$. Then at least one of the elements $A + z$ or $A + 2z$ can be represented as a commutator.
  \end{lemma}

  \begin{proof}
    If $X = \theta h$ for $\theta \in \FF_3$, then from the proof of Lemma~\ref{lemma:commutator_width} it is clear that both $A + z$ and $A + 2z$ are representable as commutators. Therefore, we will assume that $X \neq \theta h$. Then $X$ can be represented as a commutator $[Y,Z]$, where $Y, Z \in \sl_2(\FF_3)$ and the matrix $Y$ is nondegenerate. Let $B = Y + \alpha x_1 + \alpha' y_1$ and $C = Z + \beta x_1 + \beta' y_1$. The condition $[B,C] \in \lbrace A + z, A + 2z \rbrace$ is equivalent to
    \begin{equation*}
      \begin{pmatrix}
        \beta \\
        \beta'
      \end{pmatrix} = Y^{-1}
      \begin{pmatrix}
        \lambda \\
        \mu
      \end{pmatrix} + Y^{-1}Z
      \begin{pmatrix}
        \alpha \\
        \alpha'
      \end{pmatrix},
    \end{equation*}
    and $\alpha\beta' - \beta\alpha' = k \in \lbrace 1,2 \rbrace$. The matrix $Y^{-1}Z = \begin{pmatrix}\gamma_{11} & \gamma_{12} \\ \gamma_{21} & \gamma_{22} \end{pmatrix}$ is not diagonal; therefore, one of the numbers $\gamma_{12}$ or $\gamma_{21}$ is nonzero. Without loss of generality, we may assume that $\gamma_{12} \neq 0$. Therefore, for $\alpha = \beta' = 0$, we obtain a quadratic equation in $\alpha'$ whose constant term depends on $k$. The discriminant of this trinomial is a square for at least one of the values of $k$, since it depends linearly on $k$. Thus the required solution of the equation exists.
  \end{proof}

  In conclusion of this section, we present the main results on the fields mentioned in Proposition~\ref{proposition:whenCommutatorWidthIs1}. Denote by $\FFF_n$ the class of fields in which every element is representable as a sum of $n$ squares. Fields from $\FFF_1$ are called quadratically closed. Such fields include, for example, all algebraically closed fields, but this class is wider.

  \begin{lemma}\label{lemma:squares_properties}
    The following assertions hold.
    \begin{enumerate}
      \item In a finite field $K$ every element is representable as a sum of two squares.
      \item Let $S_n$ be the set of elements of an infinite field $K$ which cannot be represented as a sum of $n$ squares. If the set $S_n$ is nonempty, then it is infinite.
    \end{enumerate}
  \end{lemma}

  \begin{proof}
    \begin{enumerate}
      \item For fields of characteristic $2$, the surjectivity of the Frobenius endomorphism of finite fields implies that every element is a square. If the field has odd characteristic and contains $q$ elements, then it contains exactly $\frac{q+1}{2}$ squares. Then for any $x \in K$ the intersection of the sets
      \begin{equation*}
        \lbrace x - y^2 \mid y \in K \rbrace \cap \lbrace z^2 \mid z \in K \rbrace
      \end{equation*}
      is nonempty. Therefore every element of $K$ is equal to a sum of two squares.

      \item Let $x \neq 0$ belong to $S_n$. Then $y^2x$ also belongs to $S_n$ for every nonzero $y$ in the field. Since an infinite field contains infinitely many squares, we obtain the required result.
    \end{enumerate}
  \end{proof}

  Thus all finite fields of characteristic not $2$ belong to $\FFF_2 \setminus \FFF_1$. Obviously, $\FFF_n \subseteq \FFF_m$ for $m \geqslant n$. However, not all inclusions $\FFF_n \subseteq \FFF_{n+1}$ are strict. For this purpose, we introduce the following definitions.

  \begin{definition}
    The least natural number $s$ for which
    \begin{equation*}
      -1 = a_1^2 + \cdots + a_s^2, \quad a_i \in K,
    \end{equation*}
    is called the \textit{Stufe} $s(K)$ of the field $K$. If no such $s$ exists, then by definition we set $s(K) = +\infty$.
  \end{definition}

  \begin{definition}
    The \textit{Pythagoras number} $p(K)$ of a field $K$ is the least natural number $p$ such that every sum of squares in $K$ is a sum of $p$ squares.
  \end{definition}

  The following theorems are proved in~\cite{Rajwade1993}.

  \begin{theorem}{\cite[Theorem~2.2]{Rajwade1993}}
    For any field $K$ of finite Stufe, the number $s(K)$ is a power of two. Conversely, for any power of two there exists a field $K$ with that Stufe.
  \end{theorem}

  \begin{theorem}{\cite[Theorem~16.2]{Rajwade1993}}
    Let $K$ be a field of finite Stufe $s(K)$. Then every element of $K$ is a sum of $s(K)+1$ squares, and the inequality $s(K) \leqslant p(K) \leqslant s(K)+1$ holds.
  \end{theorem}

  Thus, for any field $K$ of finite Stufe, the pair $(s(K), p(K))$ is either $(2^n, 2^n)$ or $(2^n, 2^n+1)$. In the first case, $K \in \FFF_{2^n}$; in the second case, $K \in \FFF_{2^n+1}$. In~\cite[Theorem~5.7, p.~396]{Lam2005}, it is proved that for any such pair $(s,p)$ a suitable field $K$ exists. Therefore
  \begin{equation*}
    \FFF_1 \subsetneq \FFF_2 \subsetneq \FFF_3 \subsetneq \FFF_4 \subsetneq \FFF_5 = \FFF_7 \subsetneq \FFF_8 \subsetneq \FFF_9 = \FFF_{15} \subsetneq \cdots \subsetneq \FFF_{2^{n-1}+1} = \FFF_{2^n-1} \subsetneq \FFF_{2^n} \subsetneq \cdots.
  \end{equation*}

	\section{UA-Lie rings}\label{sec:4}
  
  \begin{definition}
    A Lie ring $\RRR$ is called a \textit{UA-Lie ring} if every commutator-preserving bijection $\phi : \RRR \rightarrow \SSS$ of Lie rings is additive.
  \end{definition}

  As a natural generalization of a UA-Lie ring, one can consider the notion of an iUA-Lie ring.

  \begin{definition}
    A Lie ring $\RRR$ is called an \textit{iUA-Lie ring} if every commutator-preserving injection $\psi : \RRR \rightarrow \SSS$ of Lie rings is additive.
  \end{definition}

  In fact, one may consider arbitrary maps instead of injective ones; however, this yields an equivalent definition, since the following statement holds.

  \begin{proposition}
    A Lie ring $\RRR$ is an iUA-Lie ring if and only if every map of Lie rings $\psi : \RRR \rightarrow \SSS$ that preserves the commutator is additive.
  \end{proposition}

  \begin{proof}
    The implication from right to left is obvious. So consider an iUA-Lie ring $\RRR$ and a map $\psi : \RRR \rightarrow \SSS$ preserving the commutator. From this map one can construct a commutator-preserving injection as follows:
    \begin{equation*}
      \psi' : \RRR \rightarrow \SSS \oplus \RRR, \qquad x \mapsto (\psi(x), x).
    \end{equation*}
    Then the map $\psi'$ is additive, and hence so is $\psi$.
  \end{proof}

  In~\cite[Theorem~C]{Sosnov2026} it is proved that the algebra $\sl_2(K)$ is an iUA-Lie ring provided that $\Char K \neq 2$. The property of being an iUA-Lie ring is very strong; namely, the following proposition holds.

  \begin{proposition}{\cite[Lemma~4.3]{Sosnov2026}}\label{proposition:negForIUA}
    Let a Lie ring $\RRR$ contain an element $y$ that cannot be represented as a commutator. Then $\RRR$ is not an iUA-Lie ring.
  \end{proposition}

  In particular, there exist infinite-dimensional simple Lie algebras that are not iUA-Lie rings~\cite{Dubouloz2021}. For finite-dimensional semisimple algebras this assertion is a conjecture formulated in~\cite{Arzhantsev2025}.

  \begin{conjecture}
    Every finite-dimensional semisimple Lie algebra is an iUA-Lie ring.
  \end{conjecture}

  It is clear that if there is a map of Lie rings $\psi : \RRR \rightarrow \SSS$ preserving the commutator, and a subring $\RRR' \subseteq \RRR$ is an iUA-Lie ring, then $\psi$ is additive on $\RRR'$.

  To each commutator-preserving map $\psi : \RRR \rightarrow \SSS$ we associate its \textit{defect} $\Delta : \RRR \times \RRR \rightarrow \SSS$, defined by the formula $\Delta(a, b) = \psi(a + b) - \psi(a) - \psi(b)$.

  In~\cite{Sosnov2026}, the following properties of the defect are proved. If $\psi : \RRR \rightarrow \SSS$ is a map preserving the commutator, then for any $a, b, c \in \RRR$ we have
  \begin{enumerate}
    \item $\Delta([a, [b, c]], [b, [c, a]]) = 0$;
    \item $[\psi(a), \Delta(b, c)] = \Delta([a, b], [a, c])$;
    \item $\Delta(a + b, c) + \Delta(a, b) = \Delta(a, b + c) + \Delta(b, c)$,
  \end{enumerate}
  and if $\psi$ is a bijection, then $\Delta(a, b) \in \psi(Z(C(a)) \cap Z(C(b)))$ for any $a, b \in \RRR$, where $Z$ is the center and $C(a)$ is the centralizer of the element $a$.

  It is also proved there that if $\phi$ is a bijection, then $\phi(Z(\RRR)) = Z(\SSS)$. This implies that for any $t_1, t_2 \in Z(\RRR)$ we have $\phi(t_1) + \phi(t_2) = \phi(t_3)$, where $t_3 \in Z(\RRR)$.

  In terms of the defect, a Lie ring is a UA-Lie ring (iUA-Lie ring) if and only if for every bijection (injection) preserving the commutator, the corresponding defect is identically zero.
  
  The defect allows one to lift additivity from subgroups. More concretely, the following lemma holds.

  \begin{lemma}\label{lemma:additivityOnSubgroups}
    Suppose that a Lie ring $\RRR$ contains additive subgroups $H_1$ and $H_2$, and that a map of Lie rings $\psi : \RRR \rightarrow \SSS$ preserving the commutator has defect $\Delta$. Then, if
    \begin{equation*}
      \Delta(H_1, H_1) = \Delta(H_2, H_2) = \Delta(H_1, H_2) = 0,
    \end{equation*}
    we have $\Delta(\langle H_1, H_2 \rangle, \langle H_1, H_2 \rangle) = 0$, where the subgroup $\langle H_1, H_2 \rangle$ is generated by $H_1$ and $H_2$.
  \end{lemma}

  \begin{proof}
    It suffices to prove that $\Delta(a_1 + b_1, a_2 + b_2) = 0$ for any $a_1, a_2 \in H_1$ and $b_1, b_2 \in H_2$, but this follows from the equality
    \begin{equation*}
      \Delta(a_1 + b_1, a_2 + b_2) = \Delta(a_1 + a_2, b_1 + b_2) + \Delta(a_1, a_2) + \Delta(b_1, b_2) - \Delta(a_1, b_1) - \Delta(a_2, b_2).
    \end{equation*}
  \end{proof}

  \begin{definition}
    We call a map of Lie rings $\phi : \RRR \rightarrow \SSS$ a \textit{central map} if it preserves the commutator and its defect $\Delta(a, b)$ lies in $\phi(Z(\RRR))$ for any $a, b \in \RRR$. A central map that is a bijection will be called a \textit{central bijection}.
  \end{definition}

  Let us determine what properties central maps possess.

  \begin{lemma}\label{lemma:imad}
    Let $\phi : \RRR \rightarrow \SSS$ be a central map with defect $\Delta$, and let $a, b \in \Im(\ad x)$ for some $x \in \RRR$. Then $\Delta(a, b) = 0$.
  \end{lemma}

  \begin{proof}
    Since $a, b \in \Im(\ad x)$, there exist $y, z \in \RRR$ such that $a = [x, y]$ and $b = [x, z]$. Then
    \begin{equation*}
      \Delta(a, b) = \Delta([x, y], [x, z]) = [\phi(x), \Delta(y, z)] = 0,
    \end{equation*}
    since $\Delta(y, z) \in \phi(Z(\RRR))$.
  \end{proof}

  \begin{proposition}\label{proposition:central_properties}
    Let $\ggg$ be a Lie algebra, $\phi : \ggg \rightarrow \SSS$ a central map, and suppose that in $\ggg$ there exists a basis $\mathcal{B}$ such that for any $x, y \in \mathcal{B}$ the element $[x, y]$ decomposes with respect to $\mathcal{B}$ in the form
    \begin{equation*}
      [x, y] = \sum_{i \in I} \beta_i e_i,
    \end{equation*}
    where $\beta_i \neq 0$ for all $i$. Suppose also that there exists $z \in \ggg$ such that $e_i \in \Im(\ad z)$ for all $i$. Then
    \begin{enumerate}
      \item if an element $a = [b, c] \in \ggg$ has a decomposition with respect to the basis $\mathcal{B}$
      \begin{equation*}
        a = \sum_{i \in I_a} \alpha_i e_i,
      \end{equation*}
      then
      \begin{equation*}
        \phi(a) = \sum_{i \in I_a} \phi(\alpha_i e_i);
      \end{equation*}
      \item if the algebra $\ggg$ is perfect and has bracket width one, then the map $\phi$ is additive.
    \end{enumerate}
  \end{proposition}

  \begin{proof}
    Let $a = [b, c]$ and suppose the elements $a, b, c$ have decompositions with respect to the basis $\mathcal{B}$
    \begin{equation*}
      a = \sum_{i \in I} \alpha_i e_i, \quad b = \sum_{i \in I} \beta_i e_i, \quad c = \sum_{i \in I} \gamma_i e_i.
    \end{equation*}
    Then
    \begin{equation*}
      \phi(a) = [\phi(b), \phi(c)] = \left[\sum_{i \in I_b} \phi(\beta_i e_i) + \phi(t_b), \sum_{i \in I_c} \phi(\gamma_i e_i) + \phi(t_c)\right],
    \end{equation*}
    where $t_b, t_c \in Z(\ggg)$, because $\phi$ is a central map. Hence,
    \begin{equation*}
      \phi(a) = \sum_{i \in I_b} \sum_{j \in I_c} \phi(\beta_i\gamma_j[e_i, e_j]).
    \end{equation*}
    The condition on $\mathcal{B}$ means that $\phi(\beta_i\gamma_j[e_i, e_j])$ decomposes into a sum of images of basis elements. Using additivity on proportional elements, we obtain assertion (1).

    Now suppose that the algebra $\ggg$ is perfect and has bracket width one. Consider the elements
    \begin{equation*}
      a = \sum_{i \in I_a} \alpha_i e_i, \quad b = \sum_{i \in I_b} \beta_i e_i.
    \end{equation*}
    Then, by the first part, it follows that
    \begin{equation*}
      \phi(a) = \sum_{i \in I_a} \phi(\alpha_i e_i), \quad \phi(b) = \sum_{i \in I_b} \phi(\beta_i e_i).
    \end{equation*}
    Adding and using additivity on proportional elements, we obtain
    \begin{equation*}
      \phi(a) + \phi(b) = \sum_{i \in I_a \setminus I_b} \phi(\alpha_i e_i) + \sum_{i \in I_a \cap I_b} \phi((\alpha_i + \beta_i)e_i) + \sum_{i \in I_b \setminus I_a} \phi(\beta_i e_i),
    \end{equation*}
    which is the decomposition of $a + b$ with respect to the basis $\mathcal{B}$. By virtue of the first part, this means that $\phi(a) + \phi(b) = \phi(a + b)$, that is, the map $\phi$ is additive.
  \end{proof}

  \begin{example}
    In the Lie algebra $\sl_n(K)$, a basis satisfying the conditions of Proposition~\ref{proposition:central_properties} is the standard basis $\mathcal{B} = \lbrace E_{ij} \mid 1 \le i \neq j \le n \rbrace \cup \lbrace E_{11} - E_{ii} \mid 2 \le i \le n \rbrace$. The commutator of any two elements of $\mathcal{B}$ is either proportional to a matrix unit or a diagonal matrix, and all diagonal matrices lie in the image of the adjoint action of the Jordan cell $J_n(0)$.
  \end{example}

  As a consequence, we obtain the following assertion.

  \begin{corollary}\label{corollary:when_central_is_additive}
    Let $\ggg$ be a perfect Lie algebra with bracket width one, in which there exists a basis $\mathcal{B}$ such that the commutator of any two elements of $\mathcal{B}$ is proportional to an element of $\mathcal{B}$. Then every central map $\phi : \ggg \rightarrow \SSS$ is additive.
  \end{corollary}

  \begin{question}
    Is it true that every central map $\phi : \ggg \rightarrow \SSS$ from a perfect Lie algebra $\ggg$ with bracket width one is additive?
  \end{question}

	\section{Proof of Theorem~\ref{theorem:when_schrodinger_UA}}\label{sec:5}

  First, let us prove the lemma on when the Schr{\"o}dinger algebra $\sch_n(K)$ is not a UA-Lie ring.

  \begin{lemma}\label{lemma:schWhenNotUA}
    Let $\Char K = 2$, or let $K \neq \FF_3$ and suppose that $K$ contains an element that cannot be represented as a sum of $n$ squares. Then the Schr{\"o}dinger algebra $\sch_n(K)$ is not a UA-Lie ring.
  \end{lemma}

  \begin{proof}
    In the case $\Char K = 2$, the algebra $\sch_n(K)$ satisfies the conditions of \cite[Proposition~1]{Arzhantsev2025}, since
    \begin{equation*}
      [\sch_n(K), \sch_n(K)] \neq \sch_n(K), Z(\sch_n(K)) = \langle z \rangle, \ [\sch_n(K), \sch_n(K)] \cap Z(\sch_n(K)) = \langle z \rangle,
    \end{equation*}
    and therefore $\sch_n(K)$ is not a UA-Lie ring.

    Now suppose that $\Char K \neq 2$ and that the element $k$ cannot be represented as a sum of $n$ squares. If the field $K$ is finite, then $n = 1$, so from the condition $K \neq \FF_3$ it follows that there exists another element $s \in K$ which also cannot be represented as a square. If, however, $n > 1$, then by item~1 of Lemma~\ref{lemma:squares_properties} the field $K$ is infinite, and hence, by item~2 of Lemma~\ref{lemma:squares_properties}, there exists another element $s \in K$ which also cannot be represented as a sum of $n$ squares.

    By Lemma~\ref{lemma:when_elem_is_not_commutator}, the elements $e - \frac{k}{2}z$ and $e - \frac{s}{2}z$ cannot be represented as commutators; therefore, the bijection of $\sch_n(K)$ interchanging $e - \frac{k}{2}z$ and $e - \frac{s}{2}z$ preserves the commutator but is not additive.
  \end{proof}

  \begin{proposition}\label{proposition:sch_additivityOnV}
    A commutator-preserving injection of Lie rings $\psi : \sch_n(K) \rightarrow \SSS$ is additive on the space $\mathcal{V} = \langle f, h, e, x_i, y_i \mid 1 \leqslant i \leqslant n \rangle$ for any field $K$ with $\Char K \neq 2$.
  \end{proposition}

  We fix the following notation. Let $\alpha, \beta \in K^n$. By the expression $\alpha \mathbf{x}$ we mean the element $\alpha_1 x_1 + \cdots + \alpha_n x_n$, and similarly for $\beta \mathbf{y}$.

  \begin{lemma}\label{lemma:sch_additivityOnXY}
    A commutator-preserving injection $\psi : \sch_n(K) \rightarrow \SSS$ is additive on the space $\langle x_i, y_i \mid 1 \leqslant i \leqslant n \rangle$.
  \end{lemma}

  \begin{proof}
    By Lemma~\ref{lemma:additivityOnSubgroups}, it suffices to prove that for any $\alpha, \beta \in K^n$ the following hold:
    \begin{enumerate}
      \item $\Delta(\alpha \mathbf{x} - \alpha \mathbf{y}, \beta \mathbf{x} - \beta \mathbf{y}) = 0$;
      \item $\Delta(\alpha \mathbf{x} - \alpha \mathbf{y}, \beta \mathbf{x} + \beta \mathbf{y}) = 0$;
      \item $\Delta(\alpha \mathbf{x} + \alpha \mathbf{y}, \beta \mathbf{x} + \beta \mathbf{y}) = 0$.
    \end{enumerate}

    From \cite[Theorem~C]{Sosnov2026} it follows that the injection $\psi$ is additive on $\sl_2(K)$. Therefore
    \begin{equation*}
      \Delta(\alpha \mathbf{x}, \beta \mathbf{y}) = [\Delta(e, f), \psi(\beta\mathbf{x} + \alpha\mathbf{y})] = 0.
    \end{equation*}

    Denote by $A, B$ and $C$ the following matrices:
    \begin{equation*}
      A = 
      \begin{pmatrix}
        1 & 1 \\
        -1 & -1
      \end{pmatrix}, \qquad
      B =
      \begin{pmatrix}
        1 & -1 \\
        -1 & -1
      \end{pmatrix}, \qquad
      C = 
      \begin{pmatrix}
        1 & -1 \\
        1 & -1
      \end{pmatrix}.
    \end{equation*}

    Then
    \begin{equation*}
      \Delta(\alpha \mathbf{x} - \alpha \mathbf{y}, \beta \mathbf{x} - \beta \mathbf{y}) = [\psi(A), \Delta(\alpha \mathbf{x}, \beta \mathbf{y})] = 0,
    \end{equation*}
    \begin{equation*}
      \Delta(\alpha \mathbf{x} - \alpha \mathbf{y}, \beta \mathbf{x} + \beta \mathbf{y}) = [\psi(B), \Delta(\alpha \mathbf{x}, -\beta \mathbf{y})] = 0,
    \end{equation*}
    \begin{equation*}
      \Delta(\alpha \mathbf{x} + \alpha \mathbf{y}, \beta \mathbf{x} + \beta \mathbf{y}) = [\psi(C), \Delta(\alpha \mathbf{x}, -\beta \mathbf{y})] = 0.
    \end{equation*}
  \end{proof}

  \begin{proof}[Proof of Proposition~\ref{proposition:sch_additivityOnV}]
    By Lemmas~\ref{lemma:additivityOnSubgroups} and \ref{lemma:sch_additivityOnXY}, it suffices to prove that for any matrix $A \in \sl_2(K)$ and $\alpha, \beta \in K^n$ we have $\Delta(A, \alpha \mathbf{x} + \beta \mathbf{y}) = 0$. Let $\mathcal{W}_1 = \langle x_i \mid 1 \leqslant i \leqslant n \rangle$, $\mathcal{W}_2 = \langle y_i \mid 1 \leqslant i \leqslant n \rangle$ and $\mathcal{W} = \mathcal{W}_1 \oplus \mathcal{W}_2$. We will prove this assertion step by step.

    \textit{Step 1.} Let us prove that $\Delta(\langle h \rangle \oplus \mathcal{W}_1 \oplus \mathcal{W}_2, \langle h \rangle \oplus \mathcal{W}_1 \oplus \mathcal{W}_2) = 0$. Let $0 \neq \gamma \in K$. Then
    \begin{equation*}
      \Delta(\gamma h + \alpha \mathbf{x}, \beta \mathbf{y}) = [\Delta(e, f), \psi(\gamma f + \beta \mathbf{x} + \alpha \mathbf{y})] = 0,
    \end{equation*}
    \begin{equation*}
      \Delta(\gamma h, \alpha \mathbf{x}) = [\Delta(-f, e), \psi(\gamma e + \alpha \mathbf{y})] = 0,
    \end{equation*}
    that is, $\Delta(\langle h \rangle \oplus \mathcal{W}_1, \mathcal{W}_2) = 0$ and $\Delta(\langle h \rangle, \mathcal{W}_1) = 0$; therefore, from Lemmas~\ref{lemma:sch_additivityOnXY} and \ref{lemma:additivityOnSubgroups} we obtain the required result.

    \textit{Step 2.} Let us prove that $\Delta(ae + bf, \mathcal{W}) = 0$ for any $0 \neq a, b \in K$. We have
    \begin{equation*}
      \Delta(ae + bf, \alpha \mathbf{x} + \beta \mathbf{y}) = \left[\psi\left(-\frac{a}{2}e + \frac{b}{2}f\right), \Delta\left(h, \frac{2\beta}{b} \mathbf{x} - \frac{2\alpha}{a} \mathbf{y}\right)\right] = 0.
    \end{equation*}

    \textit{Step 3.} Let us prove that $\Delta(ae + bf + ch, \mathcal{W}) = 0$, where $0 \neq a, b \in K$ and $c \in K$. The matrix $ae + bf + ch$ can be represented in the form
    \begin{equation*}
      \begin{pmatrix}
        c & a \\
        b & -c
      \end{pmatrix} = 
      \left[
        \begin{pmatrix}
          1 & -\frac{2c}{b} \\
          0 & -1
        \end{pmatrix},
        \begin{pmatrix}
          0 & \frac{a}{2} \\
          -\frac{b}{2} & 0
        \end{pmatrix}
      \right] = [A, B].
    \end{equation*}
    Then
    \begin{equation*}
      \Delta(ae + bf + ch, \alpha \mathbf{x} + \beta \mathbf{y}) = [\psi(A), \Delta(B, \alpha' \mathbf{x} + \beta' \mathbf{y})],
    \end{equation*}
    where $\alpha', \beta' \in K^n$ are chosen so that the identity
    \begin{equation*}
      A
      \begin{pmatrix}
        \alpha'^{\top} \\
        \beta'^{\top}
      \end{pmatrix} = 
      \begin{pmatrix}
        \alpha^{\top} \\
        \beta^{\top}
      \end{pmatrix} \Leftrightarrow [A, \alpha' \mathbf{x} + \beta' \mathbf{y}] = \alpha \mathbf{x} + \beta \mathbf{y}
    \end{equation*}
    holds.

    \textit{Step 4.} Let us prove that $\Delta(bf + ch, \alpha \mathbf{x} + \beta \mathbf{y}) = 0$ for any $0 \neq c \in K$ and $b \in K$. For this, it suffices to represent the matrix $bf + ch$ in the form $[A, B]$, where the matrix $B$ is nondegenerate and $A$ has the form from Step 3. Then
    \begin{equation*}
      \Delta(bf + ch, \alpha \mathbf{x} + \beta \mathbf{y}) = [\Delta(A, \alpha' \mathbf{x} + \beta' \mathbf{y}), \psi(B)].
    \end{equation*}
    
    \textit{Case 1.} If $b^2 + 8c \neq 0$, then
    \begin{equation*}
      A = 
      \begin{pmatrix}
        -\frac{b}{2} & c \\
        1 & \frac{b}{2}
      \end{pmatrix}, \qquad
      B = 
      \begin{pmatrix}
        -\frac{b}{2} & c \\
        2 & \frac{b}{2}
      \end{pmatrix}.
    \end{equation*}

    \textit{Case 2.} Otherwise, if $b^2 + 8c = 0$, then
    \begin{equation*}
      A = 
      \begin{pmatrix}
        -\frac{b}{2} & c \\
        -1 & \frac{b}{2}
      \end{pmatrix}, \qquad
      B = 
      \begin{pmatrix}
        -\frac{b}{2} & c \\
        0 & \frac{b}{2}
      \end{pmatrix}.
    \end{equation*}

    Similarly, one can prove that $\Delta(ae + ch, \alpha \mathbf{x} + \beta \mathbf{y}) = 0$ for any $0 \neq c \in K$ and $a \in K$.
    
    \textit{Step 5.} It remains to prove that $\Delta(\lambda f, \alpha \mathbf{x} + \beta \mathbf{y}) = 0$ and $\Delta(\lambda e, \alpha \mathbf{x} + \beta \mathbf{y}) = 0$ for any $\lambda \in K$. Let us prove the first equality; the second is proved analogously. We have
    \begin{equation*}
      \Delta(\lambda f, \alpha \mathbf{x} + \beta \mathbf{y}) = [\psi(h), \Delta(A, \alpha \mathbf{x} - \beta \mathbf{y})] = 0,
    \end{equation*}
    where the matrix $A$ has the form
    \begin{equation*}
      A = 
      \begin{pmatrix}
        1 & 0 \\
        -\frac{\lambda}{2} & -1
      \end{pmatrix}.
    \end{equation*}
  \end{proof}

  \begin{lemma}\label{lemma:sch_additivityOnZ}
    For any $\lambda_1, \lambda_2 \in K$ we have $\Delta(\lambda_1 z, \lambda_2 z) = 0$.
  \end{lemma}

  \begin{proof}
    This equality follows from the formula
    \begin{equation*}
      \Delta(\lambda_1 z, \lambda_2 z) = [\psi(\lambda_1 x_1 - \lambda_2 y_1), \Delta(y_1, x_1)] = 0.
    \end{equation*}
  \end{proof}

  Finally, we are ready to prove Theorem~\ref{theorem:when_schrodinger_UA}.

  \begin{proof}[Proof of Theorem~\ref{theorem:when_schrodinger_UA}]
    Lemma~\ref{lemma:schWhenNotUA} proves that item~(2) is true. We first prove the equivalence of items~\ref{item:1}--\ref{item:3}.

    By Proposition~\ref{proposition:whenCommutatorWidthIs1}, items~\ref{item:1} and \ref{item:2} are equivalent. Also, by Proposition~\ref{proposition:negForIUA} and the perfectness of $\sch_n(K)$ for $\Char K \neq 2$, item~\ref{item:3} implies item~\ref{item:2}. To prove the implication $\ref{item:2} \Rightarrow \ref{item:3}$, we fix a commutator-preserving injection of Lie rings $\psi : \sch_n(K) \rightarrow \SSS$. By Lemma~\ref{lemma:additivityOnSubgroups}, it suffices to prove that $\Delta(A, \lambda z) = 0$ for any $A \in \mathcal{V}$ and $\lambda \in K$. Since the bracket width of $\sch_n(K)$ is equal to one, we can write $A + \lambda z$ as
    \begin{equation*}
      A + \lambda z = [B + \alpha_1 \mathbf{x} + \beta_1 \mathbf{y}, C + \alpha_2 \mathbf{x} + \beta_2 \mathbf{y}],
    \end{equation*}
    where $B, C \in \sl_2(K)$ and $\alpha_i, \beta_i \in K^n$. Therefore, by Lemma~\ref{lemma:sch_additivityOnZ} and Proposition~\ref{proposition:sch_additivityOnV},
    \begin{multline*}
      \psi(A + \lambda z) = [\psi(B) + \psi(\alpha_1 \mathbf{x} + \beta_1 \mathbf{y}), \psi(C) + \psi(\alpha_2 \mathbf{x} + \beta_2 \mathbf{y})] = \\ = \psi([B, C] + [B, \alpha_2 \mathbf{x} + \beta_2 \mathbf{y}] + [\alpha_1 \mathbf{x} + \beta_1 \mathbf{y}, C]) + \psi(\lambda z) = \psi(A) + \psi(\lambda z).
    \end{multline*}

    Thus, the equivalence of items~\ref{item:1}--\ref{item:3} is proved. Taking Lemma~\ref{lemma:schWhenNotUA} into account, it remains only to prove that the algebra $\sch_1(\FF_3)$ is a UA-Lie ring. Let us fix a commutator-preserving bijection of Lie rings $\phi : \sch_1(\FF_3) \rightarrow \SSS$. Since Proposition~\ref{proposition:sch_additivityOnV} and Lemma~\ref{lemma:sch_additivityOnZ} hold, it suffices to prove that $\Delta(X, \lambda z) = 0$ for any $X \in \mathcal{V}$ and $\lambda \in \FF_3$.

    Suppose that there exist $X \in \mathcal{V}$ and $0 \neq \lambda \in \FF_3$ such that $\Delta(X, \lambda z) \neq 0$. Then $\Delta(X, \lambda z) = \phi(\lambda z)$, that is,
    \begin{equation*}
      \phi(X + \lambda z) = \phi(X) + \phi(\lambda z) + \phi(\lambda z) = \phi(X) + \phi(2\lambda z).
    \end{equation*}
    From Lemma~\ref{lemma:commutatorOfSch1F3} it follows that the element $X + 2\lambda z$ is representable as a commutator; therefore $\phi(X) + \phi(2\lambda z) = \phi(X + 2\lambda z)$ by Proposition~\ref{proposition:central_properties}. We obtain a contradiction with the bijectivity of $\phi$. Hence the algebra $\sch_1(\FF_3)$ is a UA-Lie ring.
  \end{proof}


\end{document}